\documentclass[12pt, a4paper]{article}
\usepackage[intlimits]{amsmath}
\usepackage{amsfonts,amsthm, amssymb}
\usepackage{fancyhdr}

\newcommand{\Rbb}{\mathbb{R}}
\newcommand{\Cbb}{\mathbb{C}}
\newcommand{\Kg}{\mathfrak{K}}
\renewcommand{\tilde}{\widetilde}

\theoremstyle{plain}
\newtheorem{thm}{Theorem}[section]

\newtheorem{lem}[thm]{Lemma}
\newtheorem{lemma}[thm]{Lemma}

\theoremstyle{definition}
\newtheorem{defn}[thm]{Definition}
\newtheorem{conj}[thm]{Conjecture}
\newtheorem{exmp}[thm]{Example}
\theoremstyle{remark}
\newtheorem{rem}[thm]{Remark}

\numberwithin{equation}{section}

\DeclareMathOperator{\const}{const}
\newcommand{\dsp}{\displaystyle}
\renewcommand{\phi}{\varphi}

\begin{document}
\title{On the principal eigenvalue of a Robin problem with a large parameter}
\author{
Michael Levitin\\
\normalsize\small Department of Mathematics, Heriot-Watt University\\
\normalsize\small Riccarton, Edinburgh EH14 4AS, U.~K.\\
\normalsize\small email {\sffamily M.Levitin@ma.hw.ac.uk}
\and
Leonid Parnovski\\
\normalsize\small Department of Mathematics, University College London\\
\normalsize\small Gower Street, London WC1E 6BT, U.~K.\\
\normalsize\small email {\sffamily Leonid@math.ucl.ac.uk}}
%
%
%
\date{\sl{{\small Dedicated to Viktor Borisovich Lidski\u{\i} on the occasion of his 80th birthday}}\\\ \\
17 January 2005}
\maketitle

\begin{abstract}
We study the asymptotic behaviour of the principal eigenvalue of a Robin (or generalised Neumann) problem with a
large parameter in the boundary condition for the Laplacian in a piecewise smooth domain. We show that the
leading asymptotic term depends only on the singularities of the boundary of the domain, and give either explicit
expressions or two-sided estimates for this term in a variety of situations.
\end{abstract}

\section{Introduction}
Let $\Omega$ be an open bounded set in $\Rbb^m$ ($m \geq 1$) with piecewise smooth, but not necessarily connected,
boundary
$\Gamma:=\partial \Omega$.
We investigate the spectral boundary value problem
\begin{equation}\label{eq1.1}
- \Delta u = \lambda u\qquad\text{in }\Omega,
\end{equation}

\begin{equation}\label{eq1.2}
\frac{\partial u}{\partial n} - \gamma Gu= 0 \qquad \text{on } \Gamma.
\end{equation}

In \eqref{eq1.1}, \eqref{eq1.2}, $\frac{\partial}{\partial n}$ denotes the outward unit normal derivative,
$\lambda$ is the spectral parameter, $\gamma$ is a positive parameter (which we later on assume to be large),
and $G:  \Gamma \rightarrow \Rbb$ is a given continuous function.  We will always assume that

\begin{equation}\label{eq1.3}
\sup_{y\in \Gamma}G(y) > 0\,.
\end{equation}

We treat the problem \eqref{eq1.1}, \eqref{eq1.2} in the variational sense, associating it with the Rayleigh quotient
\begin{equation}\label{eq1.4}
\mathcal{J} (v; \gamma, G) :=
\frac{\dsp \int_\Omega |\nabla v|^2 dx - \gamma \int_{\Gamma} G |v|^2 ds}
{ \dsp \int_\Omega |v|^2 dx}\,,\qquad v \in H^1(\Omega)\,,\ v\not\equiv 0\,.
\end{equation}

For every fixed $\gamma$, the problem \eqref{eq1.1}, \eqref{eq1.2} has a discrete spectrum of eigenvalues accumulating 
to $+\infty$.
By
\begin{equation}\label{eq:Lam}
\Lambda (\Omega; \gamma, G) := \inf_{v\in H^1(\Omega)\,,\ v\not\equiv 0} \mathcal{J}(v; \gamma, G)
\end{equation}
we denote the bottom of the spectrum of \eqref{eq1.1}, \eqref{eq1.2}.

Our aim is to study the asymptotic behaviour of $\Lambda(\Omega; \gamma, G)$ as $\gamma \to + \infty$ and its
dependence upon the singularities of the boundary $\Gamma$.

The problem {\eqref{eq1.1}}--{\eqref{eq1.2}} naturally arises in the study of reaction-diffusion
equation where a distributed absorption competes with a boundary source, see \cite{LOSS, LOS} for details.

\begin{rem} Sometimes, we shall also consider {\eqref{eq1.1}}--{\eqref{eq1.2}} for an \emph{unbounded} domain $\Omega$.
In this case, we can no longer guarantee either the discreteness of the spectrum of {\eqref{eq1.1}}--{\eqref{eq1.2}}, 
or its semi-boundedness below. We shall still use, however, the notation \eqref{eq:Lam}, allowing, in principle, 
for  $\Lambda (\Omega; \gamma, G)$ to be equal to $-\infty$.
\end{rem}

\section{Basic properties of the principal eigenvalue}

We shall mostly concentrate our attention on the case of constant boundary weight $G\equiv 1$; in this case, we 
shall denote for brevity
$$
\mathcal{J} (v; \gamma) := \mathcal{J} (v; \gamma, 1)\,,\qquad \Lambda (\Omega; \gamma):=\Lambda (\Omega; \gamma, 1)\,.
$$
See Remark \ref{rem:G} for the discussion of the case of an arbitrary smooth $G\not\equiv 1$.

We start with citing the following simple result of \cite{LOS}:

\begin{lemma} For any bounded and sufficiently smooth $\Omega\subset\Rbb^m$,
$\Lambda(\Omega; \gamma)$ is a real analytic concave decreasing function of $\gamma\ge 0$,
$\left.\Lambda\right|_{\gamma = 0} = 0$, and
$$
\left.\frac{d}{d\gamma} \Lambda(\Omega; \gamma)\right|_{\gamma = 0} = -\frac{|\Gamma|_{m-1}}{|\Omega|_m}\,.
$$
\end{lemma}

The problem \eqref{eq1.1}--\eqref{eq1.2} with $G\equiv 1$ admits a solution by separation of variables in several 
simple cases.

\begin{exmp}\label{ex:1}
For a ball $B_m(0,1)=\{|x|<1\}\subset\Rbb^m$, $\Lambda = \Lambda(B_m(0,1); \gamma)$ is given implicitly by
\begin{equation*}
\begin{split}
\sqrt{-\Lambda} \tanh \sqrt{-\Lambda} = \gamma\,,\qquad &m = 1\,,\\
\sqrt{-\Lambda} \frac{I_{m/2}(\sqrt{-\Lambda})}{I_{m/2-1}(\sqrt{-\Lambda})} = \gamma\,,\qquad &m \ge 2\,,
\end{split}
\end{equation*}
where $I$ denotes a modified Bessel function.
This implies that for any
ball $B(a, R):=\{x:|x-a|<R\}\subset\Rbb^m$,
$$
\Lambda(B(a,R); \gamma) = -\gamma^2 + O(\gamma^2)\,,\qquad\gamma \to +\infty
$$
(independently of the dimension $m$ and radius $R$); it may be shown that the same asymptotics holds for an annulus
$A_m(R_1, R) = \{ |x|\in (R_1, R)\}$.
\end{exmp}

\begin{exmp}\label{ex:2}
For a parallelepiped $P(l_1,\dots,l_m):= \left\{ |x_j| < l_j : j=1,\dots,m \right\}\subset\Rbb^m $ we get
$$
\Lambda(P(l_1,\dots,l_m); \gamma) = -\sum_{j = 1}^m \frac{\mu_j^2}{l_j^2},
$$
where $\mu_j > 0$ solves a transcendental equation
$$
\mu_j \tanh \mu_j = \gamma l_j\,.
$$
Thus we obtain
$$
\Lambda(P(l_1,\dots,l_m); \gamma) = -m\gamma^2 + O(\gamma^2)\,,\qquad\gamma \to +\infty\,.
$$
\end{exmp}

\begin{exmp}\label{ex:3} Let $\Omega = (0,+\infty)$, and $\Gamma = \{0\}$. It is easy to see that 
the bottom of the spectrum is an eigenvalue $\Lambda((0,+\infty); \gamma)=-\gamma^2$, the corresponding eigenfunction 
being $\exp(-\gamma x)$. Thus we arrive at a useful (and well-known) inequality
\begin{equation}\label{eq:use}
\int_0^\infty |v'(x)|^2\,dx-\gamma (v(0))^2\ge - \gamma^2\int_0^\infty |v(x)|^2\,dx\,,
\end{equation}
valid for all $v\in H^1((0,+\infty))$.
\end{exmp}

A slightly more complicated  example is that of a planar angle
$U_\alpha:=\{z=x+iy\in\Cbb:|\arg z|<\alpha\}$ of size $2\alpha$.

\begin{exmp}\label{ex:4} Let $\Omega = U_\alpha$ with $\alpha<\pi/2$. Again the spectrum
is not purely discrete; moreover, the separation of variable does not produce a complete set of generalised
eigenfunctions. However, one can find an eigenfunction 
$u_0(x,y)=\exp(-\gamma x/\sin \alpha)$
and compute an eigenvalue $\lambda=-\gamma^2 \sin^{-2}\alpha$ explicitly. Thus 
$\Lambda(U_\alpha; \gamma)\le-\gamma^2 \sin^{-2}\alpha$. We shall now prove that this eigenvalue is in fact the 
bottom of the spectrum.

\begin{lem}\label{lem:bot} If $\alpha<\pi/2$,
\begin{equation}\label{eq:lalp}
\Lambda(U_\alpha; \gamma)=-\gamma^2 \sin^{-2}\alpha\,.
\end{equation}
\end{lem}

\begin{proof}
It is sufficient to show that for all $v\in H^1(U_\alpha)$, we have
\begin{equation}\label{eq:est5}
\int_{U_\alpha} |\nabla v|^2 dz - \gamma \int_{\partial U_\alpha} |v|^2 ds
\ge -\gamma^2 (\sin^{-2}\alpha) \int_{U_\alpha} |v|^2 dz\,.
\end{equation}
As $ds = dy/\sin\alpha$, the left-hand side of \eqref{eq:est5} is bounded below by
$$
\int dy \left(\int \left|\frac{\partial v}{\partial x}\right|^2 dx -
\frac{\gamma}{\sin\alpha} |v|^2\right)\,.
$$
For each $y$ the integrand is not smaller than $-\gamma^2 (\sin^{-2}\alpha) \int |v|^2 dx$ by
\eqref{eq:use}. Integrating over $y$ gives \eqref{eq:est5}.
\end{proof}
\end{exmp}

\begin{exmp}\label{ex:5} Let us now consider the case of an angle $U_\alpha$ with $\alpha\in[\pi/2,\pi)$.

\begin{lem}\label{lem:bott} If $\alpha\ge\pi/2$,
\begin{equation}\label{eq:lalp1}
\Lambda(U_\alpha; \gamma)=-\gamma^2\,.
\end{equation}
\end{lem}

\begin{proof}
To prove an estimate above, we for simplicity consider a rotated angle 
$\tilde U_\alpha:=\{z=x+iy\in\Cbb:0<\arg z<2\alpha\}$.
In order to get an upper bound $\Lambda(\tilde U_\alpha; \gamma)\le -\gamma^2$, 
we construct a test function in the following manner.
Let $\psi(s)$ be a smooth nonnegative function such that 
$\psi(s)= 1$ for $|s|<1/2$, and  $\psi(s)= 0$ for $|s|>1$
Set now 
$$
\chi_\tau(s)=
\begin{cases}
1\,,&\qquad\text{if }|s|<\tau-1\,,\\
\psi(|s|-(\tau-1))\,,&\qquad\text{if }\tau-1\le |s|<\tau\,,\\
0\,,&\qquad\text{otherwise}
\end{cases}
$$ 
(a parameter $\tau$ is assumed to be greater than 1).
Consider the function 
$$
v_\tau(x,y)=e^{-\gamma y} \chi_\tau(x\gamma-\tau)\,.
$$
Then one can easily compute that
\begin{equation*}
\begin{split}
\mathcal{J} (v_\tau; \gamma) &= 
\gamma^2\left(-1+\frac{\int_{-\infty}^\infty |\chi_\tau'(s)|^2\,ds}{\int_{-\infty}^\infty |\chi_\tau(s)|^2\,ds}\right)\\
&= \gamma^2\left(-1+\frac{\int_{-1}^1 |\psi'(s)|^2\,ds}{\int_{-1}^1 |\psi(s)|^2\,ds+2(\tau-1)}\right)\,,
\end{split}
\end{equation*}
and therefore $\mathcal{J} (v_\tau; \gamma)\to -\gamma^2$ as $\tau\to\infty$. Thus, 
$\Lambda(U_\alpha; \gamma)=\Lambda(\tilde U_\alpha; \gamma)\le -\gamma^2$.

To finish the proof, we need only to show that for $v\in H^1(U_\alpha)$, 
\begin{equation}\label{eq:est6}
\int_{U_\alpha} |\nabla v|^2 dz - \gamma \int_{\partial U_\alpha} |v|^2 ds
\ge -\gamma^2 \int_{U_\alpha} |v|^2 dz\,.
\end{equation}
Denote $V_\alpha = \{z:\alpha-\pi/2 < |\arg z| <\alpha\}\subset U_\alpha$. The estimate
\eqref{eq:est6} will obviously be proved if we establish
$$
\int_{V_\alpha} |\nabla v|^2 dz - \gamma \int_{\partial U_\alpha} |v|^2 ds
\ge -\gamma^2 \int_{V_\alpha} |v|^2 dz\,.
$$
But this can be done as in the proof of Lemma~\ref{lem:bot}, by integrating first along $\partial U_\alpha$, and then 
using one-dimensional inequalities
\eqref{eq:use} in the direction orthogonal to $\partial U_\alpha$.
\end{proof}
\end{exmp}

We now consider a generalization of two previous examples to the multi-dimensional case.

\begin{exmp}\label{ex:6} Let $K\subset\Rbb^m=\{x: x/|x|\in M\}$ be a cone with the cross-section $M\subset S^{m-1}$. 
Any homothety $f: x\mapsto ax$ ($x\in\Rbb^m$, $a>0$) maps $K$ onto
itself. Then, as easily seen by a change of variables $w=\gamma^{-1}x$,
\begin{equation}\label{eq:conescale}
\Lambda(K; \gamma) = \gamma^2\Lambda(K; 1)\,.
\end{equation}

In particular, if $K$ contains a half-space, then, repeating the argument of Lemma~\ref{lem:bott} with minor 
adjustments, one can show that $\Lambda(K,1)=-1$ and so 
\begin{equation}\label{eq:lamisone}
\Lambda(K; \gamma) = -\gamma^2\,.
\end{equation}
\end{exmp}

All the above examples suggest that in general one can expect
\begin{equation}\label{eq:gen}
\Lambda(\Omega; \gamma) = -C_\Omega\gamma^2 + O(\gamma^2)\,,\qquad\gamma \to +\infty\,.
\end{equation}

Some partial progress towards establishing \eqref{eq:gen} was already achieved in \cite{LOS}.
In particular, the following Theorems were proved.
\begin{thm}\label{thm:1}
 Let $\Omega \subset \Rbb^m$ be a domain with piecewise smooth boundary $\Gamma$.  Then
$$
\limsup_{\gamma \rightarrow + \infty} \frac{\Lambda(\Omega; \gamma,  1)}{\gamma^2} \leq -1.
$$
\end{thm}

\begin{thm}\label{thm:2}
Let $\Omega \subset \Rbb^m$ be a domain with smooth boundary
$\partial \Omega$.  Then
$$
\Lambda(\Omega; \gamma) = -\gamma^2(1+o(1))\,,\qquad\gamma \to +\infty\,.
$$
\end{thm}

\begin{rem} The actual statements in \cite{LOS} are slightly weaker than the versions above, but the proofs can be 
easily
modified. Note that the proof of Theorem~\ref{thm:1} can be done by constructing a test function very similar to the 
one used in the proof of 
Lemma~\ref{lem:bott}.  
\end{rem}

The situation, however, becomes more intriguing even in dimension two, if $\Gamma$ is not smooth.
Suppose that $\Omega \subset \Rbb^2$ is a planar domain with $n$ corner points
$y_1, \dots, y_n$ on its boundary $\Gamma$.
The following conjecture was made in \cite{LOS}:

\begin{conj}\label{c:conj}
Let $\Omega \subset \Rbb^2$ be a planar domain with $n$ corner points $y_1, \dots, y_n$ on its boundary $\Gamma$ and
let $\alpha_j$, $j = 1, \dots , n$ denote the inner half-angles of the boundary at the points $y_j$. Assume that
$0<\alpha_j<\frac{\pi}{2}$.
Then \eqref{eq:gen} holds with
$$
C_\Omega =
\max_{j = 1, \dots, n} \left\{\sin^{-2} (\alpha_j)\right\}\,.
$$
\end{conj}

This conjecture was proved in \cite{LOS} only in the model case when $\Omega$ is a triangle.

As we shall see later on, formula \eqref{eq:gen} does not, in general, hold if we allow $\Gamma$ to  have zero angles (i.e.,
outward pointing cusps, see Example~\ref{ex:cusp}). We shall thus restrict ourselves to the case when $\Omega$ is piecewise smooth
in a  suitable sense, see below for the precise definition. Under this assumption, we first of all prove that the asymptotic
formula \eqref{eq:gen} holds.  Moreover, we compute $C_\Omega$ explicitly in the planar case, thus proving Conjecture~\ref{c:conj}.
In the case of dimension $m\ge 3$, we give  some upper and lower bounds on $C_\Omega$, which, in some special cases, amount to a
complete answer.

\section{Reduction to the boundary}

We shall only consider the case when $\Omega$ is piecewise smooth in the following sense: 
for each point $y\in\Gamma$ there exists an infinite 
``model'' cone $K_y$  such that for a small enough ball $B(y,r)$ of radius $r$ centred at $y$  
there exists an infinitely smooth 
diffeomorphism $f_{y}:K_y\cap B(0,r)\to\Omega\cap B(y,r)$
with $f_{y}(0)=y$ and the derivative of $f_{y}$ at $0$ being the identity matrix 
(we shall write in this case that $\Omega \sim K_y$ near a point $y\in\Gamma$). 
For example, if $y$ is a regular point of $\Gamma$, then $K_y$ is  a  half-space.

We require additionally that $\Omega$ satisfies the uniform interior cone condition \cite{GT}, i.e. there exists 
a fixed cone $K$ with non-empty interior such that each $K_y$ contains a cone congruent to $K$. 
(See Example~\ref{ex:cusp} for a discussion of a case where this condition fails.)

\begin{defn}\label{defn:Cy} 
Let $\Omega \sim K_y$ near a point $y\in\Gamma$. We denote $C_y:=-\Lambda(K_y; 1)$.
\end{defn}

Our main result indicates that the asymptotic behaviour of $\Lambda(\Omega; \gamma, 1)$ is in a sense 
``localised'' on the boundary. 

\begin{thm}\label{th:lam} Let $\Omega$ be piecewise smooth  in the above sense and satisfy the uniform
interior cone condition. Then
\begin{equation}\label{eq:mainformula}
\Lambda(\Omega; \gamma) = -\gamma^2\sup_{y\in\Gamma} C_y + o(\gamma^2)\,,\qquad\gamma \to +\infty\,.
\end{equation}
\end{thm}

\begin{rem}\label{rem:G}
This result can be easily generalised for the case of our original setting of a non-constant
boundary weight $G(y)$ satisfying  \eqref{eq1.3}:
$$
\Lambda(\Omega; \gamma, G) = -\gamma^2\sup_{\substack{y\in\Gamma\\ G(y)>0}} \{G(y)^2 C_y\} +
o(\gamma^2)\,,\qquad\gamma \to +\infty\,.
$$
\end{rem}

\begin{exmp}\label{ex:cusp} 
Formula \eqref{eq:gen}  does not, in general, hold if $\Gamma$ is allowed to have outward
pointing cusps. In particular, for a planar domain
$$
\Upsilon_p = \{(x,y)\in\Rbb^2 : x>0\,,\ |y|<x^p\}\,,\qquad p > 1
$$
one can show  that
$$
\Lambda(\Upsilon_p; \gamma) \le  -\const
\begin{cases}
\gamma^{2/(2-p)}\qquad&\text{for }1<p<2\,,\\
\gamma^N \quad\text{with any }N>0\qquad&\text{for }p\ge 2\,,
\end{cases}
$$
by choosing the test function $v = \exp\left(-\gamma x^{q_p}\right)$ with $q_p = 2-p$ for $1<p<2$ and
$q_p = 2$ for $p\ge 2$.
\end{exmp}

In order to provide the explicit asymptotic formula for $\Lambda(\Omega; \gamma)$ in the piecewise smooth case 
it remains to obtain the information on the 
dependence of the constants $C_y$ upon the local geometry of $\Gamma$ at $y$.

It is easy to do this, firstly, in the case of a regular boundary in any dimension, and, secondly, in the
two-dimensional case,  where the necessary information is already contained in Lemmas~\ref{lem:bot} and~\ref{lem:bott}.

\begin{thm}\label{th:Cy1}
Let $\Gamma$ be smooth at $y$. Then $C_y=1$.

Moreover, $C_y=1$ whenever there exists an $(m-1)$-dimensional
hyperplane $H_y$ passing through $y$ such that for small $r$, $B(y,r)\cap H_y \subset \overline{\Omega}$.
\end{thm}

\begin{thm}\label{th:Cy2}
Let $\Omega\subset\Rbb^2$ and let $y\in\Gamma$ be such that $\Omega \sim U_\alpha$ near $y$.
Then
$$
C_y=
\begin{cases}
1\,,\qquad&\text{if }\alpha\ge\pi/2\,;\\
\sin^{-2}\alpha\,,&\text{if }\alpha\le\pi/2\,.
\end{cases}
$$
\end{thm}

Theorems \ref{th:lam}, \ref{th:Cy1}, and \ref{th:Cy2} prove the validity
of Conjecture \ref{c:conj}.

In more general cases, we are only able to provide the two-sided estimates on $C_y$, and obtain the precise formulae 
only under rather restrictive additional assumptions. These results are collected in Section 5.

\section{Proof of Theorem~\ref{th:lam}} 

We proceed via a sequence of auxiliary Definitions and Lemmas. 

\begin{defn}\label{defn:kg}
Let $K\subset \Rbb^{m}$ be a cone with cross-section $M\subset S^{m-1}$, and let $r>0$. By $\Kg_{r}=\Kg_r(K)$ we 
denote the family of ``truncated'' cones $K_{r,R}$ such that
$$
K_{r,R} = \{x\in\Rbb^{m}: \theta:=x/|x|\in M\subset S^{m-1}\,,\ |x|<rR(\theta)\}\,,
$$
where $R: M\to [1,m]$ is a piecewise smooth function. Thus, for any $K_{r,R}\in\Kg_{r}$ we have
$$
K\cap B(0, r) \subset K_{r,R}  \subset K\cap B(0, mr)\,.   
$$
\end{defn}

Let $K_{r,R}\in\Kg_{r}$, and let $\sharp$ be an index assuming values $D$ or $N$ (which in turn stand for Dirichlet or
Neumann boundary  conditions). By $\Lambda^\sharp(K_{r,R}; \gamma)$ we denote the bottom of the spectrum of the
boundary value problem \eqref{eq1.1} considered in  $K_{r,R}$ with boundary conditions \eqref{eq1.2} on $\partial
K_{r,R}\cap \partial K =   \{x\in\partial K_{r,R}: x/|x|\in \partial M\}$ and with the boundary condition defined by
$\sharp$ on the rest of the boundary  $\{x: x/|x|\in M\,,\ |x|=R(\theta)\}$ (this boundary value problem is of course considered in the variational sense).

It is important to note that a simple change of variables as in Example~\ref{ex:6} leads to the re-scaling relations
\begin{equation}\label{eq:resc}
\Lambda^{\sharp}(K_{r,R}; \gamma) = \gamma^2\Lambda^{\sharp}(K_{r\gamma,R}; 1)\,.
\end{equation}
These formulae show that the bottoms of the spectra 
$\Lambda^{\sharp}(K_{r,R}; \gamma)$ 
are determined (modulo a multiplication by $\gamma^2$) by a single parameter $ \mu:=r\gamma$ via
$\Lambda^{\sharp}(K_{ \mu,R}; 1)$. It is therefore the latter which we proceed to study.  

The first Lemma gives a relation between the bottoms of the spectra for an infinite cone $K$ and its finite ``cut-offs''.

\begin{lem}\label{lem:cutoff}
Let $K_{r,R}\in\Kg_{r}(K)$ and let $ \mu=r\gamma$. Then, as $ \mu\to\infty$,
$$
\frac{\Lambda^{\sharp}(K_{r,R}; \gamma)}{\gamma^2} = \Lambda(K;1) + o(1)\,.
$$
\end{lem}

\begin{proof}[Proof of Lemma~\ref{lem:cutoff}]
By \eqref{eq:resc}, we need to prove that 
$$
\lim\limits_{ \mu\to\infty} \Lambda^{\sharp}(K_{ \mu,R}; 1)= \Lambda(K;1)\,.
$$
This can be done by considering a function $v\in H^1(K)$ and comparing the Rayleigh quotients $J(v; 1)$ with  ``truncated'' quotients $J(v\psi(\cdot/ \mu); 1)$, where $\psi$ is the same as in the proof of Lemma~\ref{lem:bott}. An easy but somewhat tedious computation which we omit shows that as $ \mu\to+\infty$, we have $J(v\psi(\cdot/ \mu); 1)\to J(v; 1)$, which finishes the proof.
\end{proof}

Let $y\in\Gamma$, and let $K_y$ be a cone with cross-section $M$ such that $\Omega\sim K_y$ near $y$. Let $r>0$ and
$K_{y,r,R}\in\Kg_{r}(K_y)$. We define $\Omega_{y,r,R}:=f_y(K_{y,r,R})$, and introduce the numbers 
$\Lambda^{\sharp}(\Omega_{y,r,R}; \gamma)$ similarly to  $\Lambda^{\sharp}(K_{y,r,R}; \gamma)$.

\begin{lem}\label{lem:omeoff}
Let $ \mu>0$ be fixed. Then
\begin{equation}\label{eq:ratio}
\lim_{r\to+0}\frac{\Lambda^{\sharp}(K_{y,r,R};  \mu/r)}{\Lambda^{\sharp}(\Omega_{y,r,R};  \mu/r)}=1
\end{equation} 
uniformly over $y\in\Gamma$.
\end{lem}

\begin{proof}[Proof of Lemma~\ref{lem:omeoff}] Let us denote by  $\tilde\Omega_{y,r,R}$ an image of $\Omega_{y,r,R}$ under the
homothety  $h_{y,r}: z\mapsto y+r^{-1}(z-y)$. Conditions imposed on the mapping $f_y$ imply that as $r\to+0$,
$\tilde\Omega_{y,r,R}\to K_{y,1,R}$ in the following sense. The volume element of $\tilde\Omega_{y,r,R}$ at a point  
$(h_{y,r}\circ f_y)(xr)$ tends to the volume element of $K_{y,1,R}$ at point $x$, and the analogous statement holds for the area 
element of the
boundary. Since $ \mu$ is fixed, this implies that the  bottoms of the spectra $\Lambda^{\sharp}(\tilde\Omega_{y,r,R};  \mu)$ (with
boundary conditions as described above) converge to $\Lambda^{\sharp}(K_{y,1,R};  \mu)$ as $r\to+0$.  Now the same re-scaling
arguments as before imply \eqref{eq:ratio}. A simple compactness argument  shows that this convergence is uniform in $y\in\Gamma$.
\end{proof}

\begin{rem}
It is easy to see that the estimates of Lemmas~\ref{lem:cutoff} and~\ref{lem:omeoff} are uniform in $R$ if we assume that all 
first
partial derivatives of $R$ are bounded by a given constant.
\end{rem}

We can now conclude the proof of Theorem~\ref{th:lam} itself. First of all, given an arbitrary positive $\epsilon$ and 
$y\in\Gamma$, we use Lemma~\ref{lem:cutoff} to find a positive $\mu(y)$ such that
\begin{equation}\label{eq:mu}
|\gamma^{-2}\Lambda^{\sharp}(K_{y,r,R}; \gamma)+C_y|<\frac{\epsilon}{2}\,,
\end{equation}
whenever $\gamma \ge r^{-1} \mu(y)$.  It is easy to see that $\mu(y)$ can be chosen to be continuous on each smooth component of the boundary. Therefore, there exists $\tilde\mu=\sup\limits_{y\in\Gamma}\mu(y)$.
Let us fix this value of $\tilde\mu$ for the rest of the proof.   

Formula \eqref{eq:mainformula} splits into two 
asymptotic inequalities. The inequality 
$$
\Lambda(\Omega; \gamma) \le -\gamma^2\sup_{y\in\Gamma} C_y + \epsilon\gamma^2\,,
\qquad\gamma \to +\infty
$$
follows immediately from formula~\eqref{eq:mu}, Lemma~\ref{lem:omeoff} (with $\sharp=D$ and $\mu=\tilde\mu$) and the obvious inequality 
$$
\Lambda(\Omega; \gamma) \le \Lambda^D(\Omega_{y,r,R}; \gamma)\,.
$$

In order to prove the opposite inequality 
\begin{equation}\label{eq:estbelow}
\Lambda(\Omega; \gamma) \ge
-\gamma^2\sup_{y\in\Gamma} C_y + \epsilon\gamma^2\,,\qquad\gamma \to +\infty\,,
\end{equation} 
we consider a partition
$\overline\Omega=\overline{\bigsqcup_{\ell=0}^N Q_\ell}$ by disjoint sets $Q_\ell$ satisfying the 
following properties: $Q_0\Subset \Omega$ (i.e. $Q_0\cap\Gamma=\emptyset$), and for each $\ell\ge 1$,
$Q_\ell=\Omega_{y,r,R}=f_y(K_{r,R})$ with some $r>0$, $y\in\Gamma$, and $K_{r,R}\in\Kg_{r}(K_y)$, such that $\Omega\sim
K_y$ near $y$. Such a partition can be constructed for each sufficiently small $r>0$ by considering, for example, a
partition of $\Rbb^m$ into cubes  of size $r$, and including into $Q_0$ all the cubes which lie strictly inside
$\Omega$. Note that  $\Gamma = \overline{\bigcup_{\ell=1}^N (\Gamma\cap Q_\ell)}$.

Now we use the following inequality: assuming that 
$J(v;\gamma)$ is negative for some $v\in H^1(\Omega)\setminus\{0\}$, we have
\begin{equation}\label{eq:estQ0}
\begin{split}
J(v;\gamma) &= \frac{\dsp \int_{\Omega} |\nabla v|^2 dx - \gamma \int_{\Gamma}  |v|^2 ds}
{ \dsp \int_\Omega |v|^2 dx} \ge
\frac{\dsp \int_{\Omega\setminus Q_0} |\nabla v|^2 dx - \gamma \int_{\Gamma}  |v|^2 ds}
{ \dsp \int_{\Omega\setminus Q_0} |v|^2 dx}\\
&=  \frac{\dsp \sum_{\ell=1}^N \int_{Q_\ell} |\nabla v|^2 dx - \gamma \sum_{\ell=1}^N 
\int_{\Gamma\cap Q_\ell}  |v|^2 ds}
{ \sum_{\ell=1}^N \dsp \int_{Q_\ell} |v|^2 dx}\\
&\ge
\min_{\ell=1\dots N} \frac{\dsp \int_{Q_\ell} |\nabla v|^2 dx - \gamma \int_{\Gamma\cap Q_\ell}  |v|^2 ds}
{ \dsp \int_{Q_\ell} |v|^2 dx}\,.
\end{split}
\end{equation}

Note that the last expression in \eqref{eq:estQ0} is bounded below by $\inf \Lambda^N(\Omega_{y,r,R}; \gamma)$, where
the infimum is taken over all $y\in\Gamma$ and all functions $R$ admissible in the sense of Definition~\ref{defn:kg}.

Finally, taking the size of the partition $r\to+0$, and using formula~\eqref{eq:mu} and Lemma~\ref{lem:cutoff} with $\mu=\tilde\mu$, we obtain \eqref{eq:estbelow}. 

\section{Estimates in the general case}

Let us now discuss the general case. As we have already shown, the problem of computing the constant
$\dsp C_\Omega=\sup_{y\in\Gamma} C_y$ in  \eqref{eq:gen} is reduced to calculating the bottoms of the spectra
$\Lambda(K_y;1)=-C_y$ for infinite model cones $K_y$. We have also shown that $C_y=1$ when $\Gamma$ is smooth at $y$.
We now consider a case when $\Gamma$ is singular  at $y$. 

Let $j$ be the co-dimension of a singularity of $\Gamma$ at $y$. By this we mean that $K_y = \Rbb^{m-j}\times \tilde{K}$,
with $\tilde{K}= \{z\in\Rbb^j: z/|z|\in\tilde{M}\}$, with the singular cross-section $\tilde{M}\subset S^{j-1}$.
If $j\ge 3$, we restrict our analysis to the case when the closure of $\tilde{M}$ is contained in open hemisphere 
$\{\theta\in S^{j-1}:\theta_1>0\}$. For simplicity, we assume that $\tilde{M}$ is convex (this stronger
requirement may be relaxed, see Remark~\ref{rem:nonconv}). 

The case $j=1$ corresponds to a regular point $y\in\Gamma$. The case $j=2$ is treated in exactly the same way as the
planar case, as in this situation  $\tilde{K}=U_\alpha$ and the constant $C_y$ is the same as in Theorem~\ref{th:Cy2}.    

Consider now the case $j\ge 3$.
It might seem natural to introduce the spherical coordinates on $\tilde{K}$ at this stage. Unfortunately, such an
approach is not likely to succeed --- although the variables separate, the resulting lower-dimensional problems
are coupled in a complicated way. Indeed, Example~\ref{ex:4} shows that the principal eigenfunction is not easily
expressed in spherical coordinates. Therefore, we will try to choose a coordinate frame more suitable for this
problem. Once more, Example~\ref{ex:4} gives us a helpful insight into what this coordinate frame should be.

We need more notation. Let $w\in\tilde{K}$ with  $\theta=w/|w|\in \tilde{M}$. We define $\Pi_{\theta}$ as
a $(j-1)$-dimensional hyperplane  passing through $\theta$ and orthogonal to $w$. Let 
$P_\theta = \Pi_\theta \cap \partial \tilde{K}$. We need to consider only the points $\theta$ such that $P_\theta$ 
is bounded and $\theta\in P_\theta$. Such directions $\theta$ always exist due to the convexity of  
$\tilde{M}$.

We now introduce the coordinates $(\xi,\eta)\in\Rbb\times\Rbb^{j-1}$ of a point $z\in\Rbb^j$, such that 
$\xi = z\cdot\theta$ is a coordinate along $\theta$  and  $\eta = z-\xi\theta$ represent coordinates along the plane
$\Pi_{\theta}$.   

We also need the spherical coordinates $(\rho, \phi)$ with the origin at $\theta$ on $\Pi_\theta$, such that
$\rho=|\eta|$ and $\phi=\eta/|\eta|\in S^{j-2}$. We define a function $b(\phi)=b_\theta(\phi)$ in such a way that
$P_\theta = \{(\rho, \phi) : \rho=b(\phi)\}$.

In these coordinates, 
\begin{equation}\label{eq:Kcoords}
\tilde{K} = \{(\xi,\rho,\phi): \xi>0\,,\ \rho<\xi b_\theta(\phi)\}
\end{equation}
and
\begin{equation}\label{eq:bKcoords}
\partial\tilde{K} = \{(\xi,\rho,\phi): \xi>0\,,\ \rho=\xi b_\theta(\phi)\}\,.
\end{equation}

Denote 
\begin{equation}\label{eq:sig3}
\sigma_\theta(\phi):=\sqrt{1+b_\theta^{-2}(\phi)+(b'_\theta(\phi))^2 b_\theta^{-4}(\phi)}\,.
\end{equation}

We are ready now to formulate a general statement in the case $j=3$.

\begin{thm}\label{th:Cy3}
Let $y\in\Gamma$ be a singular point of co-dimension three in the above sense. Then the
constant $C_y$ satisfies the following two-sided estimates:
\begin{equation}\label{eq:est}
\sup_\theta \left(\frac{\dsp\int_{S^1} b_\theta^2(\phi)
\sigma_\theta(\phi)\,d\phi}
{\dsp\int_{S^1} b_\theta^2(\phi)\,d\phi}\right)^2 \le C_y
\le
\inf_\theta\sup_\phi \sigma^2_\theta(\phi)
\end{equation}
\end{thm}

\begin{rem}\label{rem:nonconv} 
Theorem~\ref{th:Cy3} can be extended to the case of non-convex $\tilde{M}$. Then,
the function $b_\theta(\phi)$ (which defines the boundary) may become multivalued. In that case we need to treat the 
integrals
in the left-hand side of \eqref{eq:est} separately along each branch of $b_\theta$, and count them with a plus or 
minus
sign.
\end{rem}

\begin{proof}[Proof of Theorem~\ref{th:Cy3}]
The separation of variables shows that $C_y=-\Lambda(\tilde{K};1)$.
We start by estimating $C_y$ below (and thus $\Lambda(\tilde{K};1)$ above). 
Let us fix $\theta\in\tilde{M}$ satisfying the above conditions; for brevity we shall omit the subscript $\theta$ 
in all the intermediate calculations. 

Consider the following test function
\begin{equation}\label{eq:psiest}
v(z) = \exp(- a \xi)\,,\qquad z=(\xi,\rho,\phi)\in \tilde{K}\,,
\end{equation}
where $a$ is a positive parameter to be chosen later.

Then we explicitly calculate
\begin{equation}\label{eq:int3}
\int_{\tilde{K}} v^2(z)\,dz  
=\int_0^\infty \exp(- 2a\xi)\,d\xi\int_{S^1}d\phi\int_0^{\xi b(\phi)}\rho d\rho = \frac{1}{8a^3}\int_{S^1} b^2(\phi)\,d\phi
\end{equation}
and
\begin{equation}\label{eq:int1}
\int_{\tilde{K}} |\nabla v(z)|^2\,dz  
=a^2\int_{\tilde{K}} v^2(z)\,dz = \frac{1}{8a}\int_{S^1} b^2(\phi)\,d\phi\,.
\end{equation}

Let us now calculate the integral along the boundary $\partial\tilde{K}$. For each
$\tilde{\eta}=(\tilde{\rho},\tilde{\phi})\in\Rbb^2$  there exists a unique point 
$z=(\tilde{\xi},\tilde{\eta})=(\tilde{\xi}(\tilde{\eta}),\tilde{\eta})\in\partial\tilde{K}$, where one can easily
compute $\dsp\tilde{\xi}(\tilde{\eta})=\frac{\tilde{\rho}}{b(\tilde{\phi})}$. Thus the area element of the 
boundary  $ds$ can be expressed as $\dsp\frac{1}{\cos\beta}d\tilde{\eta}$, where $\beta$ is an angle between two
planes. One of these planes is  $\Pi_\theta$ and the other one is the plane containing the origin and the straight line
$L$ which lies in $\Pi_\theta$ and is tangent to  $P_\theta$ at the point $\xi=1$, $\rho=b(\tilde{\phi})$,
$\phi=\tilde{\phi}$. Without loss of generality we assume now that $\tilde{\phi}=0$, otherwise we just rotate the
picture. Then the equation of $L$ in cartesian coordinates $\eta=(\eta_1,\eta_2)$ on $P_\theta$ becomes 
$L=\{\eta_1=b(0)+tb'(0)\,,\ \eta_2=b(0)t: t\in\Rbb\}$. It is a simple geometric exercise to show that the base of the
perpendicular  dropped from the origin onto $L$ corresponds to the parameter value $\dsp
t^*=-\frac{b(0)b'(0)}{b(0)^2+(b'(0))^2}$ and therefore this base point is given by  $\dsp (\eta^*_1,\eta_2^*)=
\frac{b(0)^2}{b(0)^2+(b'(0))^2}(b(0),-b'(0))$. Another geometric exersise shows that $\cot\beta$ is equal to the length
of the vector $(\eta^*_1,\eta_2^*)$, and therefore
$$
\frac{1}{\cos\beta}=\sqrt{1+\cot^{-2}\beta}=\sqrt{1+b^{-2}(0)+(b'(0))^2 b^{-4}(0)}\,.
$$
Thus, the area element, with account of \eqref{eq:sig3}, is
\begin{equation}\label{eq:ds}
ds = \frac{1}{\cos\beta}d\tilde{\eta} = 
\sqrt{1+b^{-2}(\tilde{\phi})+(b'(\tilde{\phi}))^2 b^{-4}(\tilde{\phi})}\,d\tilde{\eta}=
\sigma(\tilde{\phi})\,d\tilde{\eta}\,,
\end{equation}
and we can evaluate the boundary contribution as
\begin{equation}\label{eq:int2}
\begin{split}
\int_{\partial\tilde{K}} v^2(z)\,ds &= 
\int_{\Rbb^2} \exp(- 2a\tilde{\xi})\sigma(\tilde{\phi})\,d\tilde{\eta}\\
&= \int_{S^1}\,d\tilde{\phi}\,\sigma(\tilde{\phi})
\int_{0}^\infty \tilde{\rho}\exp(- 2a\tilde{\rho}/b(\tilde{\phi}))\,d\tilde{\rho}\\
&= \frac{1}{4a^2}\int_{S^1}b^2(\phi)\sigma(\phi)\,d\phi 
\,.
\end{split}
\end{equation} 

Combining now \eqref{eq:int3}, \eqref{eq:int1}, and \eqref{eq:int2}, we obtain
$$
J(v;1) = a^2 - \frac{\dsp 2a\int_{S^1}b^2(\phi)\sigma(\phi)\,d\phi}{\dsp\int_{S^1} b^2(\phi)\,d\phi}\,.
$$

Optimising with respect to $a$ gives 
\begin{equation}\label{eq:a3}
a=\frac{\dsp\int_{S^1} b_\theta^2(\phi)
\sigma(\phi)\,d\phi}
{\dsp\int_{S^1} b_\theta^2(\phi)\,d\phi}\,,
\end{equation}
and further optimization with respect to $\theta$ produces the desired lower bound in \eqref{eq:est}.

Let us now prove the upper bound on $C_y$ in \eqref{eq:est}, which corresponds to the lower bound on 
$\Lambda(\tilde{K};1)$.
We need to show that for any $v\in H^1(\tilde{K})$ and any $\theta\in\tilde{M}$ the following inequality holds:
\begin{equation}\label{eq:a4}
\int_{\tilde{K}}|\nabla v(z)|^2\,dz-\int_{\partial\tilde{K}}|v(z)|^2\,ds\ge 
-\left(\sup_\phi\sigma(\phi)\right)^2\int_{\tilde{K}}|v(z)|^2\,dz\,.
\end{equation}
Using the obvious estimate
$$
\int_{\tilde{K}} |\nabla v|^2\,dz \ge \int_{\tilde{K}} |\partial_{\xi} v|^2\,dz\,,
$$
formula \eqref{eq:ds} for the area element, and inequality \eqref{eq:use} in the variable $\xi$ for 
each value of $\eta$, we arrive at 
\eqref{eq:a4}. This completes the proof.
\end{proof}	

\begin{rem} In the case of a three-edged corner (i.e. when $\tilde{M}$ is a two-dimensional spherical triangle) the
left- and right-hand sides of \eqref{eq:est} in fact coincide, so Theorem \ref{th:Cy3} gives the exact expression for
$C_y$. The same is true if $\tilde{M}$ is a spherical polygon which has an inscribed circle (i.e., a circle touching
all the sides of $\tilde{M}$). Indeed, in this case the supremum in the left-hand side and the infimum in the
right-hand side of \eqref{eq:est} are equal and are attained when $\theta$ is the centre of the inscribed circle. This
immediately follows from the fact that in this case and for this choice of $\theta$,  $\sigma\equiv\text{\rm const}$.
Moreover, it is easy to see that the test function \eqref{eq:psiest} with the parameter $a$ given by \eqref{eq:a3} is
an eigenfunction with the eigenvalue at the bottom of the spectrum $\Lambda(\tilde{K}; 1)$.

Thus, Theorems \ref{th:Cy1},~\ref{th:Cy2}, and~\ref{th:Cy3} provide an exact asymptotics
of $\Lambda(\Omega; \gamma)$ whenever $m=3$ and each vertex of $\Omega$ has three edges coming from it.
\end{rem}

Assume now that $j>3$. This case is pretty much similar to the previous one (in particular, the test function used in 
obtaining the estimate below on $C_y$ is still given by \eqref{eq:psiest}), the only difference being that the  area
element of the boundary now becomes a volume element of co-dimension one and is much more cumbersome to calculate. 
We skip the detailed calculations.

In order to state the result, we need more notation.
Define a $(j-2)$-dimensional vector
$\zeta_\theta(\phi) := b_\theta(\phi)\nabla_\phi b_\theta(\phi)$ and  $(j-2)\times(j-2)$
matrix $Z_\theta(\phi) := b^2_\theta(\phi) I + (\nabla_\phi\otimes\nabla_\phi) b_\theta(\phi)$. Now put
$\Psi_\theta(\phi):= (Z^{-1}_\theta(\phi)\zeta_\theta(\phi))$ and
\begin{equation}\label{eq:Sig}
\Sigma_\theta(\phi) := \sqrt{1+\left((b_\theta(\phi)-\Psi_\theta(\phi)\cdot\nabla_\phi b_\theta(\phi))^2+
b^2_\theta(\phi)|\Psi_\theta(\phi)|^2\right)^{-1}}\,.
\end{equation}

\begin{thm}\label{th:Cy4}
Let $y\in\Gamma$ be a singular point of co-dimension $j\ge 4$ in the above sense. Then the
constant $C_y$ satisfies the following two-sided estimates:
\begin{equation}\label{eq:est4}
\sup_\theta \left(\frac{\dsp\int_{S^{j-2}} b_\theta^{j-1}(\phi)
\Sigma_\theta(\phi)\,d\phi}
{\dsp\int_{S^{j-2}} b_\theta^{j-1}(\phi)\,d\phi}\right)^2 \le C_y
\le
\inf_\theta\sup_\phi\Sigma^2_\theta(\phi)\,.
\end{equation}
\end{thm}

\begin{rem} It is easily seen that Theorem~\ref{th:Cy3} is in fact a partial case of  Theorem~\ref{th:Cy4} if we 
formally set $j=3$ in the latter.
Indeed, for $j=3$ all the quantities depend upon a scalar parameter $\phi$, and we obtain
$$
\zeta_\theta(\phi) = b_\theta(\phi) b'_\theta(\phi)\,,\quad
Z_\theta(\phi) = b^2_\theta(\phi) + (b'_\theta(\phi))^2\,,\quad
\Psi_\theta(\phi) = \frac{b_\theta(\phi) b'_\theta(\phi)}{b^2_\theta(\phi) + (b'_\theta(\phi))^2}\,,
$$
giving
\begin{align*}
\Sigma_\theta(\phi) &=
\sqrt{1+
\left(\left(b_\theta(\phi)-\frac{b_\theta(\phi) (b'_\theta(\phi))^2}{b^2_\theta(\phi) + (b'_\theta(\phi))^2}\right)^2
+b^2_\theta(\phi)\frac{(b'_\theta(\phi))^2}{(b^2_\theta(\phi) + (b'_\theta(\phi))^2)^2}\right)^{-1}}\\
&=\sqrt{1+\frac{b^2_\theta(\phi) + (b'_\theta(\phi))^2}{b^4_\theta(\phi)}}=\sigma_\theta(\phi)\,,
\end{align*}
so that formula~\eqref{eq:est4} becomes \eqref{eq:est}.
\end{rem}

\begin{rem} As before, the estimates \eqref{eq:est4} give the precise value of $C_y$ whenever $M$ is a
$(j-1)$-dimensional spherical polyhedron which admits an inscribed ball (for example when $M$ has exactly $j$
faces). Moreover, the bottom of the spectrum is again an eigenvalue corresponding to the eigenfunction 
\eqref{eq:psiest}.
\end{rem}

\

{\bf Acknowledgements:} our collaboration was partially supported by the EPSRC Spectral Theory
Network.

\

\end{document}